\documentclass[11pt]{elsarticle}

\usepackage[english]{babel}
\usepackage{csquotes}
\usepackage{hyperref}
\usepackage{nicematrix}

\usepackage[T1]{fontenc}

\usepackage{amssymb}
\usepackage{dcolumn}
\usepackage{enumerate}
\usepackage{mathdots}
\usepackage{tabularx}
\usepackage{subfigure}
\usepackage{amsmath}
\usepackage{amsthm}
\usepackage{tikz} 
\usepackage{float}
\usepackage{mathtools}
\usepackage{graphicx}
\usepackage{pgfplots}
\usepackage{natbib}

\usepackage[margin=3.5cm]{geometry}

\usepackage{pgf,tikz}
\usepackage{mathrsfs}

\usepackage{algorithm,algpseudocode,float}
\allowdisplaybreaks

\newcommand{\pmn}{\mathbb{P}_m(n)}

\newcommand{\sign}{\text{sign}}

\newcommand{\indsize}{\scriptsize}
\newcommand{\colind}[2]{\displaystyle\smash{\mathop{#1}^{\raisebox{.5\normalbaselineskip}{\indsize #2}}}}
\newcommand{\rowind}[1]{\mbox{\indsize #1}}

\newcommand{\vl}{\vec{v}_\lambda}

\usepackage{pgf,tikz}\usepackage{mathrsfs}\usetikzlibrary{arrows}

\usetikzlibrary{matrix,fit, arrows, shapes}

\numberwithin{equation}{section}

\newtheorem{thm}{Theorem}[section]
\newtheorem{prop}[thm]{Proposition}

\newtheorem{lem}[thm]{Lemma}

\newtheorem{conj}[thm]{Conjecture}

\theoremstyle{definition}
\newtheorem{defn}{Definition}[section]

\pgfplotsset{compat=1.18}

\begin{document}
 
\title{On the spectra of prefix-reversal graphs}

\author[1]{Sa\'ul A. Blanco}
\affiliation[1]{organization={Department of Computer Science, Indiana University}, city={Bloomington, IN},
postcode={47408}, country={U.S.A.}}
\ead{sblancor@iu.edu}

\author[2]{Charles Buehrle}
\affiliation[2]{organization={Department of Mathematics, Physics, and Computer Studies, Notre Dame of Maryland University}, city={Baltimore, MD}, postcode={21210}, county={U.S.A.}}
\ead{cbuehrle@ndm.edu}

\date{June 9, 2025}

\begin{abstract} 
In this paper, we study spectral properties of prefix-reversal graphs. These graphs are obtained by connecting two elements of $C_m\wr S_n$ via prefix reversals. If $m=1,2$, the corresponding prefix-reversal graphs are the classic pancake and burnt pancake graphs. If $m>2$, then one can consider the directed and undirected versions of these graphs. We prove that the spectrum of the undirected prefix-reversal graph $\pmn$ contains all even integers in the interval $[0,2n]\setminus\{2\lfloor n/2\rfloor\}$ and if $m\equiv0\pmod4$, we then show that the spectrum contains all even integers in $[0,2n]$. In the directed case, we show that the spectrum of the directed prefix-reversal graph $P(m,n)$ contains all integers in the interval $[0,n]\setminus\{\lfloor n/2\rfloor\}$. As a consequence, we show that in either case, the prefix-reversal graphs have a small spectral gap. 

\end{abstract}

\maketitle

\section{Introduction}

The \emph{pancake problem} is a classic sorting problem first proposed in 1975 by Jacob Goodman, (under the pen name Harry Dweighter --``harried waiter'')~\cite{Dweighter75} where one is tasked with sorting a stack of pancakes, all of which are of a different diameter, utilizing a spatula. To perform the task, one is only allowed to use the spatula to pick up an arbitrary number of pancakes from the top of the stack, flip them, and place them back into the stack. The number of pancakes that may be flipped can vary each time. Gates and Papadimitriou provided the first non-trivial algorithm to solve it~\cite{GatesPapa} as well as proposing the \emph{burnt pancake problem} where, in addition, each pancake in the stack has an orientation (burnt side) and the goal is to sort a stack in such a way that at the end each pancake is in the right order according to size \textit{and} orientation (burnt side down). 

While finding the \emph{optimal} number of flips to sort any stack of pancakes is computationally challenging, Bulteau, Fertin, and Rusu showed it to be NP-hard~\cite{BulFerRusu}, the time complexity to sort an arbitrary stack of pancakes with spatula flips is $\Theta(n)$. 

Given that in the pancake problem each pancake has a different diameter, it is natural to represent each stack of $n$ pancakes with a permutation of length $n$. In this context, flipping pancakes corresponds to performing a \emph{prefix reversal} on a permutation written in one-line notation. A natural object of study is the so-called \emph{pancake graph}, where two permutations are connected by an edge if and only if there is a prefix reversal that transforms one into the other. One can define the \emph{burnt pancake graph} for the burnt pancake problem if one utilizes \emph{signed permutations} instead of permutations. Several structural properties of these graphs are known, in particular, regarding their cycle structure, some references include the authors along with Patidir, Kanevsky and Feng, as well as Konstantinova and Medvedev~\cite{BB23, BBP19, BBP19Perm, KF95,KM16}.

Formally, the group of signed permutations can be thought of as $C_2\wr S_n$, the wreath product of the order-two group and the symmetric group. This interpretation motivates considering an extension of the pancake and burnt pancake graphs using elements of $C_m\wr S_n$, with $m>2$, as vertices and connecting two elements whose permutations are a prefix reversal apart and the reversed entries are incremented in their respective ``colors,'' cycling modulo $m$. We provide the full formal definitions in Section~\ref{sec:pmn}. These graphs are referred to as \emph{prefix-reversal graphs}. If $m=1,2$ the corresponding prefix-reversal graphs are the pancake graph $P_n$ and the burnt pancake graph $BP_n$, respectively. 

Regarding their spectra, Dalf\'o and Fiol~\cite{Dalfo} established that the spectrum of $P_n$ contains the set $[-1,n-1]\setminus\{\lfloor n/2\rfloor -1\}$, where the notation is taken to be the integer interval. Moreover, the authors proved that the spectrum of $BP_n$ contains the set $[0,n]\setminus\{\lfloor n/2\rfloor\}$. In this article, we establish a similar result for prefix-reversal graphs with $m>2$.

\subsection{Main results} We extend the authors' main result~\cite{BB24} and establish the existence of a range of integer eigenvalues in the prefix-reversal graphs, both when they are considered directed and undirected. More precisely, we prove that 
\begin{enumerate}
    \item The spectrum of the directed prefix-reversal graph, denoted by $P(m,n)$, contains the set $[0,n]\setminus\{\lfloor n/2\rfloor\}$.
    \item The undirected prefix-reversal graph, denoted by $\pmn$, has a spectrum that contains every even integer in the set $[0,2n]\setminus\{2\lfloor n/2\rfloor\}$. Furthermore, in the case $m\equiv 0\pmod 4$, the spectrum of $\pmn$ contains every even integer in $[0,2n]$.
\end{enumerate}

The general behavior of multiplicities of these eigenvalues for prefix-reversal graphs remains unknown. 

\subsection{Organization of the paper}

In Section~\ref{sec:prelim}, we include all the needed formal definitions and notations used throughout the paper. We discuss notions of prefix-reversal graphs, regular partitions vertices and their connection to graph spectra, and circulant matrices, and determinants of block matrices when the blocks are circulant matrices. 

In Section~\ref{sec:permutationmatrices} we describe the permutation matrices corresponding to prefix reversals, which we use in Section~\ref{sec:undirected} and Section~\ref{sec:directed} where we present our results for the undirected case $\pmn$ and the directed case $P(m,n)$, respectively. Finally, in Section~\ref{sec:final}, we conclude with some final remarks, including a few conjectures. 

\section{Preliminaries}\label{sec:prelim}

\subsection{Basic notation}

If $a,b$ are two integers with $a<b$, then we denote the set $\{a,a+1,\ldots, b\}$ by $[a,b]$, and if $a=1$, then we simply use $[b]$. For a positive integer $n$, we let $S_n$ denote the symmetric group of degree $n$. Similarly, for a positive integer $m$, we let $C_m$ denote the cyclic group of order $m$ and write its elements additively $\{0,1,\ldots, m-1\}$. Lastly, we will use $\mathbf{0}_m$ to denote the $m\times m$ matrix that has a 0 in each of its entries, and we will write $I_m$ to denote the $m\times m$ identity matrix. 

\subsection{Circulant matrices}  Let us further recall that a \emph{circulant matrix} is a matrix where each row is a cyclic shift (one position to the right) of the previous row. In other words, if $C$ is an $\ell\times\ell$ circulant matrix, then $C$ has the following form:
\[
C =
\begin{bmatrix} 
c_0 & c_1 & c_2 & \cdots & c_{\ell-1} \\ 
c_{\ell-1} & c_0 & c_1 & \cdots & c_{\ell-2} \\ 
c_{\ell-2} & c_{\ell-1} & c_0 & \cdots & c_{\ell-3} \\ 
\vdots & \vdots & \vdots & \ddots & \vdots \\ 
c_1 & c_2 & c_3 & \cdots & c_0  
\end{bmatrix}.
\]

For convenience in notation, we will denote the circulant matrix above by $C(c_0,c_1,\ldots,c_{\ell-1})$. The following result is well-known (see, for example, Varga~\cite{Varga54}). 
\begin{lem}\label{lem:circulant-eigenvalues}
    The eigenvalues of the circulant matrix $C(c_0,c_1,\ldots,c_{\ell-1})$ have the form 
    \[
    \lambda_{k}=\sum_{j=0}^{\ell-1}c_ie^{-2\pi ikj/\ell},
    \]
    where $i$ is the complex unit and $k \in [0, \ell-1]$.
\end{lem}

Furthermore, it is well known that the set of $\ell \times \ell$ circulant matrices, $\mathcal{C}^\ell$, with matrix addition and multiplication, forms a commutative subring of the ring of all $\ell \times \ell$-matrices with complex entries, $\mathbb{C}^{\ell \times \ell}$; see, for example, I. Kra and S. Santiago~\cite{IKraSSant}. This fact allows one to find the determinant of a block matrix whose blocks are circulant matrices, by the following theorem, see J. Silvester~\cite[Theorem 1]{silvester:hal-01509379}.
\begin{thm}\label{thm:Silvester}
    If $M$ is an $mn\times mn$ matrix containing $n^2$ blocks of $m\times m$ circulant matrices, then \[\det(M) = \det\left(\det_{\mathcal{C}^m}(M)\right).\] 
\end{thm}
That is, the determinant of the matrix can be found by taking the determinant of the matrix that results from treating the blocks of the matrix as entries. As a simple example, consider the block matrix
    \[
        M = \begin{bmatrix}
            A & B \\ C & D
        \end{bmatrix},
    \]
    where $A, B, C,$ and $D$ are each $m \times m$ circulant matrices, then the determinant of $M$ is 
    \[
        \det(M) =\det(AD-BC),
    \]
    where $AD-BC$ is an $m \times m$ circulant matrix.

\subsection{Prefix-reversal graphs}\label{sec:pmn}

Let us write $S(m,n)$ to refer to the group $C_m\wr S_n$ of \emph{colored permutations}. The elements of $C_m$ are referred to as \emph{colors}. The entries in the elements of $S_n$ are referred to as \emph{symbols}. Each symbol $\pi_j$ of the permutation $\pi \in S_n$ has an associated color $e_j \in C_m$. We use superscript notation to pair symbol and color. Thus, $\pi \in S(m,n)$ can be written as $\pi = \pi_1^{e_1}\pi_2^{e_2} \cdots\pi_n^{e_n}$.

Given $\pi=\pi^{e_1}_1\pi^{e_2}_2\cdots\pi_n^{e_{n}}\in S(m,n)$, we define the $i^{\text{th}}$-\emph{flip} of $\pi$, denoted by $r_i^+(\pi)$ as
\[
r^+_i(\pi^{e_1}_1\pi^{e_2}_2\cdots\pi_n^{e_{n}})=\pi^{a_i}_i\pi^{a_{i-1}}_{i-1}\cdots\pi_1^{a_1}\pi_{i+1}^{e_{i+1}}\cdots\pi_n^{e_{n}},
\] where $a_j=e_j+1$ mod $m$, for $1\leq j\leq i$.

Moreover, we define the $i^{\text{th}}$-\emph{flop} of $\pi$, denoted by $r_i^-(\pi)$, as
\[
r^-_i(\pi^{e_1}_1\pi^{e_2}_2\cdots\pi_n^{e_{n}})=\pi^{b_i}_i\pi^{b_{i-1}}_{i-1}\cdots\pi_1^{b_1}\pi_{i+1}^{e_{i+1}}\cdots\pi_n^{e_{n}},
\] where $b_j=e_j-1$ mod $m$, for $1\leq j\leq i$.

It is also useful to have an alternative notation in mind. Following Steingr\'imsson~\cite{S94} (who uses the term \emph{indexed permutations} to refer to the elements of $S(m,n)$), we denote an element of $S(m,n)$ as a product $s\times \sigma$ where $s$ is an $n$-tuple from the $n$-th fold product of $C_m$ with itself, and $\sigma\in S_n$. For example, $3^15^02^04^21^0\in S(3,5)$ can also be denoted by $(1,0,0,2,0)\times 35241$. The product notation makes it easy to understand how multiplication in $S(m,n)$ works. Indeed, if $s\times\sigma_1,r\times\sigma_2\in S(m,n)$, their product $(s\times\sigma_1)\cdot (r\times\sigma_2)$ is given by $(s+\sigma_1(r))\times \sigma_1(\sigma_2)$, where $\sigma_1(r)=(r_{\sigma_2(1)},r_{\sigma_1(2)},\ldots,r_{\sigma_1(n)})$ and the addition is done mod $m$. 

With the product notation of $S(m,n)$ as well as disjoint cycle notation for the symmetric group, the $i^{\text{th}}$-flip considered as an element in $S(m,n)$ is  
\[
r^+_i = (\underbrace{1,1,\ldots,1}_i,\underbrace{0,0,\ldots,0}_{n-i})\times \left((1\;i)(2\;i-1)\cdots\left(\left\lfloor\frac{i+1}{2}\right\rfloor\;\left\lceil\frac{i+1}{2}\right\rceil\right)\right),
\] where $(a\;b)$, with $1\leq a<b\leq n$, denotes the transposition that swaps $a$ and $b$ and leaves every other symbol fixed. Similarly, the $i^{\text{th}}$-flop is
\[
r^-_i = (\underbrace{n-1,n-1,\ldots,n-1}_i,\underbrace{0,0,\ldots,0}_{n-i})\times \left((1\;i)(2\;i-1)\cdots\left(\left\lfloor\frac{i+1}{2}\right\rfloor\;\left\lceil\frac{i+1}{2}\right\rceil\right)\right),
\]

Therefore, we can also think of either type of prefix reversals as closely related colored permutations. For example, $r_1^+(1^j)=1^{j+1\text{ mod }m}$, and in general, for $j\in [i]$ and $k\in C_m$,
\[
r^+_i(j^k)=\left((1\;i)(2\;i-1)\cdots\left(\left\lfloor\frac{i+1}{2}\right\rfloor\;\left\lceil\frac{i+1}{2}\right\rceil\right)\right)(j)^{k+1\text{ mod }m} = (i-j+1)^{k+1\text{ mod }m},
\] otherwise when $j>i$ we have $r^+_i(j^k)=j^k$. That is, $r^+_i$ maps $j\in[i]$ to the product of transpositions $(1\;i)(2\;i-1)\cdots\left(\left\lfloor\frac{i+1}{2}\right\rfloor\;\left\lceil\frac{i+1}{2}\right\rceil\right)$ applied to $j$, and adds 1 to $k$ mod $m$. Similarly,
\[
r^-_i(j^k)=\left((1\;i)(2\;i-1)\cdots\left(\left\lfloor\frac{i+1}{2}\right\rfloor\;\left\lceil\frac{i+1}{2}\right\rceil\right)\right)(j)^{k-1\text{ mod }m} = (i-j+1)^{k-1\text{ mod }m},
\] for $j \in [i]$ or $r^-_i(j^k) = j^k$ when $j>i$.

Let us denote the set of all prefix-reversals (flips and flops) by $R^\pm_n=\{r_i^+,r^-_i\mid i\in[n]\}$ and the set of all flips by $R^+_n=\{r^+_i\mid i\in [n]\}$. Then we define the \emph{prefix-reversal graphs} as follows.

\begin{defn}  
\begin{enumerate} Let $m,n$ be two positive integers. Then,
\item The undirected \emph{prefix-reversal} graph, denoted by $\pmn$, has the elements of $S(m,n)$ as vertices and the set $\{(v,rv)\mid v\in S(m,n),\, r\in R^\pm_n\}$ as the set of edges
    \item The directed \emph{prefix-reversal graph}, denoted by $P(m,n)$, has the elements of $S(m,n)$ as vertices and the set $\{(v,rv)\mid v\in S(m,n),\,r\in R^+_n\}$ as set of edges.
\end{enumerate} Notice that these prefix-reversal graphs are Cayley graphs (defined below) of $S(m,n)$ utilizing $R^\pm_n$ or $R^+_n$ as generators, respectively. Moreover, $P(1,n)\cong \mathbb{P}_1(n)\cong P_n$, the pancake graph, and $P(2,n)\cong \mathbb{P}_2(n)\cong BP_n$, the burnt pancake graph. 
\end{defn}

\subsection{Regular partitions and connections to graph spectrum} 

Recall that given a group $Gr$ with generating set $S$, the \emph{Cayley graph} of $Gr$ with respect to $S$, denoted by $Cay(Gr,S)$, is the graph $G=(V,E)$ with $V=Gr$ and $E=\{(g,sg):g\in Gr,\,s\in S\}$. When the set $S$ is closed under inverses, the resulting Cayley graph can be considered undirected. Otherwise, the Cayley graph is a directed graph.

The \emph{spectrum} of a graph $G$ refers to the set of eigenvalues of the adjacency matrix $A_G$ of $G$, along with their multiplicities. 

For a given $\pi\in S_n$, its corresponding \emph{permutation matrix} is given by $P(\pi)=(p_{i,j})_{n\times n}$ where $p_{i,j}=1$ if $j=\pi(i)$ and $p_{i,j}=0$, otherwise. We now describe a technique from Dalf\'o and Fiol~\cite{Dalfo} based on a theorem by Godsil~\cite{Godsil}. In essence, the technique can be summarized in the following proposition.

\begin{prop}(See~\cite[Lem 1.1 \& Prop 2.1]{Dalfo})\label{prop:method}
    Let $Gr$ is a subgroup of the symmetric group $S_n$ with generating set $S=\{\pi_1,\ldots,\pi_\ell\}$. Then the spectrum of the matrix $\sum_{i=1}^\ell P(\pi_i)$ is contained in the spectrum of the adjacency matrix of $Cay(Gr,S)$.
\end{prop}

For completeness, we outline the needed technical details to justify the proposition above. Let us recall that a partition $\mathcal{P}$ of a set $S$ is a collection of sets $\mathcal{P}=\{P_1,\ldots,P_k\}$ whose union is the entirety of $S$ and if $i\neq j$ then $P_i\cap P_j=\emptyset$. After selecting an ordering of the elements of $V$, each $P_i$ has a characteristic column vector $\mathbf{v}_{i}$ whose $u^{\text{th}}$ component is 1 if $u\in P_i$ and 0 if $u\not\in P_i$. The \emph{characteristic matrix} is the $|V|\times k$ matrix whose $i^{\text{th}}$ column is the vector $\mathbf{v}_{i}$. 

Given a graph $G=(V,E)$, a partition of its vertex set $V$, $\mathcal{P}=\{P_1, \ldots, P_k\}$, is called \emph{equitable} or \emph{regular} if for any $v^i_1,v^i_2\in P_i$, one has that $|P_j\cap N(v^i_1)|=|P_j\cap N(v^i_2)|$ for all $1\leq i,j\leq k$. Here, $N(v)$ denotes the \emph{neighbors} of $v$: the set $\{u\in V:(u,v)\in E\}$ of all vertices in $V$ connected to $v$ by an edge. Given an equitable partition $\mathcal{P}=\{P_1,\ldots,P_k\}$, notice that the quantities $|P_i\cap N(u)|$ does not depend on $u$ itself, depending, rather, on the partition $\mathcal{P}$. By defining $b_{i,j}=|P_j\cap N(u)|$ for $u\in P_i$, one defines the \emph{quotient matrix} of $G$ with respect to the equitable partition $\mathcal{P}$ by $(b_{i,j})_{k\times k}$.

It is not hard to prove that in the set-up of Proposition~\ref{prop:method}, there is an equitable partition $\mathcal{P}$ of the vertices of $Cay(Gr,S)$, with $S=\{\pi_1,\ldots,\pi_{\ell}\}$, where  $\sum_{i=1}^\ell P(\pi_i)$ is the quotient matrix of $Cay(Gr,S)$ with respect to $\mathcal{P}$. Indeed, one can construct $\mathcal{P}=\{P_1,\ldots,P_n\}$ where $P_i=\{\pi\in Gr: \pi(i)=a\}$, where $a$ is any element of $[n]$. It is not difficult to see $\mathcal{P}$ is an equitable partition and the matrix $\sum_{i=1}^{\ell} P(\pi_i)$ is the quotient matrix of $Cay(Gr,S)$ with respect to $\mathcal{P}$.

The last piece needed to establish Proposition~\ref{prop:method} is the following theorem that can be found in Godsil's book~\cite{Godsil}.

\begin{thm}[\cite{Godsil}] Let $G=(V,E)$ be a graph with adjacency matrix $A_G$ and let $\mathcal{P}=\{P_1,\ldots,P_\ell\}$ be a partition of $V$ with characteristic matrix $M_{\mathcal{P}}$. Then,
\begin{enumerate}
    \item[(i)] The partition $\mathcal{P}$ is regular if and only if there exits a matrix $C$ such that $M_{\mathcal{P}}C=M_{\mathcal{P}}A_G$. Furthermore, if such a $C$ exists, then $C$ is the quotient matrix of $A_G$ with respect to $\mathcal{P}$.
    \item[(ii)] If $\mathcal{P}$ is regular and $\mathbf{v}$ is an eigenvector of the quotient matrix $Q$ of $A_G$ with respect to $\mathcal{P}$, then $M_{\mathcal{P}}\mathbf{v}$ is an eigenvector of $A_G$. In particular, then the spectrum of $B$ is contained in the spectrum of $A_G$.
\end{enumerate}
\end{thm}
    
As a direct consequence, the spectrum of $\sum_{i=1}^\ell P(\pi_i)$ is contained in the spectrum of $Cay(Gr,S)$.

Notice that we can regard $S(m,n)$ as a subgroup of $S_{mn}$ generated by prefix-reversals. Therefore, we are able to apply Proposition~\ref{prop:method} with $Gr=S(m,n)$ and $S=R^{\pm}_n$.


\section{Permutation matrices corresponding to $r^{\pm}_{i}$}\label{sec:permutationmatrices}

Let us index the permutation matrices of prefix-reversals in $S(m,n)$ by the following ordering of the entries in a colored permutation: $i_1^{j_1}<i_2^{j_2}$ if and only if $i_1<i_2$, or $j_1<j_2$ if $i_1=i_2$. For example, the order of the entries in the elements of $S(3,2)$ is given by $1^0<1^1<1^2<2^0<2^1<2^2$. 

Given $\pi\in S(m,n)$, we denote its corresponding $mn\times mn$ permutation matrix by $P(\pi)$. The product notation allows us to easily describe the permutation matrices of prefix reversals. Indeed, $r_i^+$ will permute the symbols according to $$\left((1\;i)(2\;i-1)\cdots\left(\left\lfloor\frac{i+1}{2}\right\rfloor\;\left\lceil\frac{i+1}{2}\right\rceil\right)\right)$$ and the colors will increase by 1 mod $m$. In the case of $r_i^-$, the colors decrease by 1 mod $m$. For example, 

\[
  P(r^+_{2}) =
  \begin{array}{@{}c@{}}
    \rowind{$1^0$} \\ \rowind{$1^1$} \\ \rowind{$1^2$} \\ \rowind{$2^0$} \\ \rowind{$2^1$} \\\rowind{$2^2$}
  \end{array}
  \left[
  \begin{array}{ *{6}{c} }
     \colind{0}{$1^0$}  &  \colind{0}{$1^1$}  &  \colind{0}{$1^2$}  & \colind{0}{$2^0$} & \colind{1}{$2^1$} & \colind{0}{$2^2$} \\
0 & 0 & 0  & 0 & 0 & 1 \\
 0 & 0 & 0  & 1 & 0 & 0 \\
 0 & 1  & 0  & 0 & 0 & 0 \\
 0 & 0 & 1  & 0 & 0 & 0 \\
 1 & 0 & 0  & 0 & 0 & 0
  \end{array}
  \right] \text{ and }
  P(r^-_{2}) =
  \begin{array}{@{}c@{}}
    \rowind{$1^0$} \\ \rowind{$1^1$} \\ \rowind{$1^2$} \\ \rowind{$2^0$} \\ \rowind{$2^1$} \\\rowind{$2^2$}
  \end{array}
  \left[
  \begin{array}{ *{6}{c} }
     \colind{0}{$1^0$}  &  \colind{0}{$1^1$}  &  \colind{0}{$1^2$}  & \colind{0}{$2^0$} & \colind{0}{$2^1$} & \colind{1}{$2^2$} \\
0 & 0 & 0  & 1 & 0 & 0 \\
 0 & 0 & 0  & 0 & 1 & 0 \\
 0 & 0  & 1  & 0 & 0 & 0 \\
 1 & 0 & 0  & 0 & 0 & 0 \\
 0 & 1 & 0  & 0 & 0 & 0
  \end{array}
  \right].
\]

If $B=(b_{i,j})_{k\times k}$ is a $k\times k$-matrix and $1\leq i_1<i_2\leq k$, $1\leq j_1<j_2\leq k$, then $B[(i_1,i_2),(j_1,j_2)]$ denotes the submatrix of $B$ given by $(b_{i,j})$ with $i_1\leq i\leq i_2$ and $j_1\leq j\leq j_2$.

Notice that the blocks $P(r^+_{2})[(1^0,1^2),(2^0,2^2)]$ and $P(r^+_{2})[(2^0,2^2),(1^0,1^2)]$ are the circulant matrix $C(0,1,0)$. Moreover, $P(r^-_{2})[(1^0,1^2),(2^0,2^2)]$ and $P(r^-_{2})[(2^0,2^2),(1^0,1^2)]$ are the circulant matrix $C(0,0,1)$.

These observations discussed above allow us to prove the following theorem. We write $\mbox{diag}(d_1,d_2,\ldots,d_\ell)$ to refer to the $\ell\times\ell$ diagonal matrix with diagonal entries $d_1,\ldots, d_\ell$.

\begin{thm}
    For $m>2,n>1$, let
    \begin{itemize}
\item $D(m,n)=
\mbox{diag}(\underbrace{0,\ldots,0}_m,\underbrace{2,\ldots,2}_m,\ldots,\underbrace{2(n-1),\ldots,2(n-1)}_m)$,
\item $C^\pm(m)=C(0,1,\underbrace{0,\ldots,0}_{m-2})+C(\underbrace{0,\ldots,0}_{m-1},1)$,
\item $C^+(m)= C(0,1,\underbrace{0,\ldots,0}_{m-2})$,
\item $M^\pm(m,n)=\sum_{i=1}^nP(r_i^+)+P(r_i^-)$, \text{and}
\item $M^+(m,n)=\sum_{i=1}^nP(r_i^+)$.
\end{itemize}
Furthermore, let
\[C^{\pm}(m,n)=
\begin{bmatrix}
C^\pm(m)      & C^\pm(m)  & C^\pm(m)      & \cdots      & C^\pm(m) & C^\pm(m)      \\
C^\pm(m)      & C^\pm(m)       & C^\pm(m)& \cdots      & C^\pm(m) & \mathbf{0}_{m}      \\
C^\pm(m)      & C^\pm(m)      & C^\pm(m)      & \cdots & \mathbf{0}_{m} & \mathbf{0}_{m}      \\
\vdots & \vdots & \vdots & \vdots & \ddots & \vdots \\
C^\pm(m)      & C^\pm(m)       & \mathbf{0}_{m}      & \mathbf{0}_{m}      & \mathbf{0}_{m}      & \mathbf{0}_{m} \\
C^\pm(m)      & \mathbf{0}_{m}      & \mathbf{0}_{m}      & \mathbf{0}_{m}      & \mathbf{0}_{m}      & \mathbf{0}_{m}
\end{bmatrix}_{mn\times mn},
\] and
\[C^+(m,n)=
\begin{bmatrix}
C^+(m)     & C^+(m)   & C^+(m)       & \cdots      & C^+(m)  & C^+(m)       \\
C^+(m)       & C^+(m)        & C^+(m) & \cdots      & C^+(m)  & \mathbf{0}_{m}      \\
C^+(m)       & C^+(m)       & C^+(m)       & \cdots & \mathbf{0}_{m} & \mathbf{0}_{m}      \\
\vdots & \vdots & \vdots & \vdots & \ddots & \vdots \\
C^+(m)      & C^+(m)        & \mathbf{0}_{m}      & \mathbf{0}_{m}      & \mathbf{0}_{m}      & \mathbf{0}_{m} \\
C^+(m)       & \mathbf{0}_{m}      & \mathbf{0}_{m}      & \mathbf{0}_{m}      & \mathbf{0}_{m}      & \mathbf{0}_{m}
\end{bmatrix}_{mn\times mn},
\]  where $\mathbf{0}_m$ denotes the $m\times m$ zero matrix. Then,
\begin{align*}
    M^\pm(m,n)&=C^\pm(m,n)+D(m,n) \text{ and}\\
    M^+(m,n)&=C^+(m,n)+\frac{1}{2}D(m,n)
\end{align*}
\end{thm}

In light of Theorem~\ref{thm:Silvester}, we can think of $C^\pm(m,n)$ and $C^+(m,n)$ as $mn\times mn$ block matrices formed by $m\times m$ matrices $C^\pm(m), C^+(m)$, and $\mathbf{0}_{m}$, respectively. 
    
\begin{proof}

Recall that if $B=(b_{i,j})_{k\times k}$ is a $k\times k$-matrix and $1\leq i_1<i_2\leq k$, $1\leq j_1<j_2\leq k$, then $B[(i_1,i_2),(j_1,j_2)]$ denotes the submatrix of $B$ given by $(b_{i,j})$ with $i_1\leq i\leq i_2$ and $j_1\leq j\leq j_2$. We extend the notation to our matrices indexed by $i^k$, with $1\leq i\leq n$ and $0\leq k\leq m-1$, linearly ordered as discussed above.

We will argue by cases depending on the positions of the entries of $M^\pm(m,n)$.

   \textbf{Diagonal entries of $M^\pm(m,n)$.} Since $r^+_i,r^-_i$ leave fixed every entry in a colored permutation, $k^j$, with $k>i$, then there are no $r^+_i,r^-_i$ that leave fixed entries at the diagonal positions of $M^\pm(m,n)[(1^0,1^{m-1}),(1^0,1^{m-1})]$, only $r^+_1,r^-_1$ leave the diagonal entries of $M^\pm(m,n)[(2^0,2^{m-1}),(2^0,2^{m-1})]$ fixed, and in general only $r^+_i,r^-_i$ with $i\leq k+1$ and $k\geq0$ leave the diagonal entries of $M^\pm(m,n)[(k^0,k^{m-1}),(k^0,k^{m-1})]$ fixed. Thus the diagonal entries of $M^\pm(m,n)$ are given by $D(m,n)$.
    

\textbf{Non-diagonal entries of $M^\pm(m,n)$.} Let us notice that the action of $r^+_i$ is to add 1 mod $m$ to the colors of any symbol at position $\leq i$ and that the action of $r_i^-$ to the colors of symbols at positions $\leq i$ is to subtract 1 mod $m$. 

Notice that, if $r_i^+,r^-_i$ swap position $a$ and position $b$, with $1\leq a,b\leq i$, then $a+b=i+1$ and therefore $P(r_i^+)+P(r_i^-)$ will only contribute non-zero entries to the non-diagonal entries of the blocks of the form $M^\pm(m,n)[(a^0,a^{m-1}),(b^0,b^{m-1})]$ with $a+b=i+1$. More specifically, $$M^\pm(m,n)[(a^0,a^{m-1}),(b^0,b^{m-1})]=C(0,1,\underbrace{0,\ldots,0}_{m-2})+C(\underbrace{0,\ldots,0}_{m-1},1)=C^\pm(m),$$ with $a+b=i+1$.

Moreover, the prefix-reversals $\{r^+_i,r^-_i:i\in[n]\}$ would not affect non-diagonal entries of blocks of the form $M^\pm(m,n)[(a^0,a^{m-1}),(b^0,b^{m-1})]$ with $a+b>n+1$, and therefore if $a+b>n+1$, $M^\pm(m,n)[(a^0,a^{m-1}),(b^0,b^{m-1})]=\mathbf{0}_m$.
 

Therefore,
\[
M^\pm(m,n)=C^\pm(m,n)+D(m,n). 
\] The proof for $M^+(m,n)$ is similar except that we only need to account for the $n$ flips $r^+_i$ with $1\leq i\leq n$.  Therefore, in this case the diagonal entries should be half of $D(m,n)$ and blocks of the form $C^\pm(m,n)$ are of the form $C^+(m,n)$.
\end{proof}

For example, $M^\pm(4,3)$ is given by the matrix below. We use $r^\pm_i$ as shorthand for ``$r^+_i$ and $r^-_i$''.
\[
  \begin{array}{@{}c@{}}
     \\ \rowind{$r_1^{\pm}$ changes the color of positions $[(1^0,1^3),(1^0,1^3)]$ only} \\ \rowind{$r_2^{\pm}$ changes the color of positions $[(1^0,1^3),(2^0,2^3)]$ only} \\  \rowind{$r_3^{\pm}$ changes the color of positions $[(1^0,1^3),(3^0,3^3)]$ only} 
     \\ \\ \rowind{$r_1^{\pm}$ does not change the colors of positions $[(2^0,2^3),(j^0,j^3)]$ since $2>1$} \\ \rowind{$r_2^{\pm}$ changes the color of positions $[(2^0,2^3),(1^0,1^3)]$ only} \\  \rowind{$r_3^{\pm}$ changes the color of positions $[(2^0,2^3),(2^0,2^3)]$ only}
     \\ \\ \rowind{$r_1^{\pm}$ does not change the color of positions $[(3^0,3^3),(j^0,j^3)]$ since $3>1$} \\ \rowind{$r_2^{\pm}$ does not change the color of positions $[(3^0,3^3),(j^0,j^3)]$ since $3>2$} \\\ \rowind{$r_3^{\pm}$ changes the color of positions $[(3^0,3^3),(1^0,1^3)]$ only}
  \end{array}
\left[
\begin{array}{cccc|cccc|cccc}
 0 & 1 & 0 & 1 & 0 & 1 & 0 & 1 & 0 & 1 & 0 & 1 \\
 1 & 0 & 1 & 0 & 1 & 0 & 1 & 0 & 1 & 0 & 1 & 0 \\
 0 & 1 & 0 & 1 & 0 & 1 & 0 & 1 & 0 & 1 & 0 & 1 \\
 1 & 0 & 1 & 0 & 1 & 0 & 1 & 0 & 1 & 0 & 1 & 0 \\
 \hline 
 0 & 1 & 0 & 1 & 2 & 1 & 0 & 1 & 0 & 0 & 0 & 0 \\
 1 & 0 & 1 & 0 & 1 & 2 & 1 & 0 & 0 & 0 & 0 & 0 \\
 0 & 1 & 0 & 1 & 0 & 1 & 2 & 1 & 0 & 0 & 0 & 0 \\
 1 & 0 & 1 & 0 & 1 & 0 & 1 & 2 & 0 & 0 & 0 & 0 \\
 \hline
 0 & 1 & 0 & 1 & 0 & 0 & 0 & 0 & 4 & 0 & 0 & 0 \\
 1 & 0 & 1 & 0 & 0 & 0 & 0 & 0 & 0 & 4 & 0 & 0 \\
 0 & 1 & 0 & 1 & 0 & 0 & 0 & 0 & 0 & 0 & 4 & 0 \\
 1 & 0 & 1 & 0 & 0 & 0 & 0 & 0 & 0 & 0 & 0 & 4 \\
\end{array}
\right].
\]


The following Lemma will be useful when we find the eigenvalues of $M^\pm(m,n)$ and $M^+(m,n)$ and reveals a curiosity that occurs when the number of colors is a multiple of 4.

\begin{lem}\label{lem:multipleof4}
    The matrix $C^\pm(m)$ is singular if and only if $m\equiv0\pmod4$.
\end{lem}
\begin{proof} Both $C^\pm(m)$ and $C^+(m)$ are circulant matrices. Therefore, by Lemma~\ref{lem:circulant-eigenvalues}
    the eigenvalues of $C^\pm(m)$, $\lambda_{\ell}$, have the following form
    \[
    \lambda_{\ell}= e^{\frac{-2\pi i \ell}{m}}+ e^{-\frac{2\pi i (m-1)\ell}{m}}=2 \cos\left(\frac{2\pi\ell}{m}\right).
    \] Then $\lambda_\ell=0$ for some $\ell$ if and only if $m=4k$ and $\ell=k$ for some positive integer $k$. 
\end{proof} Incidentally, the matrix $C^+(m)$ is never singular for any $m>3$.

\section{The undirected case $\pmn$}\label{sec:undirected}

We first consider the case of $Cay(S(m,n),R^\pm(n))$. In this case, the graph obtained may be considered as undirected since the generating set is $R^\pm(n)$ symmetric, that is, each generator and its inverse are in the generating set. 

\begin{thm}\label{thm:undirected}
        If $m>2$, then the spectrum of $\pmn$ contains all the even integers from the set $[0,2n]\setminus \{2\lfloor n/2\rfloor\}$. Moreover, if $m\equiv0\pmod 4$, then  the spectrum of $\pmn$ contains all the even integers from the set $[0,2n]$. 
\end{thm}

\begin{proof}  
The proof here is achieved with direct verification that the following are eigenvalue-eigenvector pairs for $M^\pm(m,n)$.

Since $\pmn$ is $2n$-regular, the eigenvalue $2n$ has an eigenvector that is the $mn$-dimensional vector with all entries equal to 1. 

    The eigenvalue $2(n-1)$ has an eigenvector of
    \[
    (\underbrace{0,\ldots,0}_{m\text{ times}},\underbrace{-1,\ldots,-1}_{m(n-2)\text {times}},\underbrace{n-2,\ldots,n-2}_{m\text{ times}})^\top,
    \] the eigenvalue $2(n-2)$ has an eigenvector of
    \[
    (\underbrace{0,\ldots,0}_{2m\text{ times}},\underbrace{-1,\ldots,-1}_{m(n-4)\text {times}},\underbrace{n-4,\ldots,n-4}_{m\text{ times}},\underbrace{0,\ldots,0}_{m\text{ times}})^\top.
    \] In general, the eigenvalue $2(n-i)$ with $1\leq i<\lfloor n/2\rfloor$ has an eigenvector of 
    \[(\underbrace{0,\ldots,0}_{mi \text{ times}},\underbrace{-1,\ldots,-1}_{m(n-2i)\text{ times}}\,\underbrace{n-2i,\ldots,n-2i}_{m\text{ times}},\underbrace{0,\ldots,0}_{m(i-1)\text{ times}})^\top. 
    \] Moreover, if $n$ is even, say $n=2\ell$ for some positive integer $\ell$, then the eigenvalue $2(\ell-i)$ with $1\leq i\leq  \ell$ has eigenvector 
    \[
    (\underbrace{0,\ldots,0}_{m(\ell-i) \text{ times}},\underbrace{-2i+1,\ldots,-2i+1}_{m\text{ times}},\underbrace{1,\ldots,1}_{m(2i-1)\text{ times}},\underbrace{0,\ldots,0}_{m(\ell-i)\text{ times}})^\top.
    \] Furthermore, if $n$ is odd, say $n=2\ell+1$ for some positive integer $\ell$, then the eigenvalue $2(\ell+1)$ has eigenvector 
    \[
    (\underbrace{0,\ldots,0}_{m\ell \text{ times}},\underbrace{-1,\ldots,-1}_{m\text{ times}},\underbrace{1,\ldots,1}_{m\text{ times}},\underbrace{0,\ldots,0}_{m(\ell-1)\text{ times}})^\top, 
    \] and finally if $n=2\ell+1$, then the eigenvalue $2(\ell+1-i)$ with $1\leq i\leq \ell+1$ has eigenvector 
    \[
    (\underbrace{0,\ldots,0}_{m(\ell-i) \text{ times}},\underbrace{-2i,\ldots,-2i}_{m\text{ times}},\underbrace{1,\ldots,1}_{2mi\text{ times}},\underbrace{0,\ldots,0}_{m(\ell-i)\text{ times}})^\top.
    \] 
    Therefore, the spectrum of $\pmn$, which contains the spectrum of $M^{\pm}(m,n)$ by Proposition~\ref{prop:method}, contains every even integer from the set $[0,2n]\setminus\{2\lfloor n/2\rfloor\}$.

    As an illustrative example, we verify that $\lambda = 2(n-i)$, for $1 \leq i \leq \lfloor n/2 \rfloor$, is an eigenvalue of $M^\pm(m,n)$ with corresponding eigenvector 
    
    \begin{align*}
        \vl &= (\underbrace{0,\ldots,0}_{mi \text{ times}},\underbrace{-1,\ldots,-1}_{m(n-2i)\text{ times}}\,\underbrace{n-2i,\ldots,n-2i}_{m\text{ times}},\underbrace{0,\ldots,0}_{m(i-1)\text{ times}})^\top\\
        &= (\underbrace{\overbrace{0, \ldots, 0}^{m\text{ times}}, \ldots, \overbrace{0, \ldots, 0}^{m\text{ times}}}_{i \text{ times}}, \underbrace{\overbrace{-1, \ldots, -1}^{m\text{ times}},\ldots,\overbrace{-1, \ldots, -1}^{m\text{ times}}}_{(n-2i)\text{ times}}, \overbrace{n-2i, \ldots, n-2i}^{m \text{ times}}, \underbrace{0,\ldots, 0}_{m(i-1)\text{ times}})^\top\\
        &= (\underbrace{\vec{0}_m, \ldots, \vec{0}_m}_{i\text{ times}}, \underbrace{\overrightarrow{(-1)}_m, \ldots, \overrightarrow{(-1)}_m}_{(n-2i)\text{ times}}, \overrightarrow{(n-2i)}_m, \underbrace{\vec{0}_m, \ldots, \vec{0}_m}_{(i-1)\text{ times}})^\top,
    \end{align*}
    where $\vec{k}_m$ is an $m$-dimensional vector with all $k$ entries.
    
    We will show that $M^\pm(m,n) \vl=\lambda \vl$. Consider $M^\pm(m,n) \vl = C^\pm(m,n) \vl + D(m,n) \vl$, and we argue by cases depending on the entries of $M^\pm(m,n) \vl $ by considering an arbitrary row of the $M^\pm(m,n)$, say the $j^{\text{th}}$ one, thought of as a block matrix. 

        \[\begin{bmatrix}
            \vdots \qquad \vdots & \vdots \qquad \vdots\\
            \underbrace{C^\pm(m) \; \cdots \; C^\pm(m)}_{n-j+1\text{ times}} & \underbrace{\mathbf{0}_m \; \cdots \; \mathbf{0}_m}_{j-1\text{ times}}\\
            \vdots \qquad \vdots & \vdots \qquad \vdots\\
        \end{bmatrix} \vl + \begin{bmatrix}
            \vdots \qquad \vdots & \vdots & \vdots \qquad \vdots\\
            \underbrace{\mathbf{0}_m \; \cdots \; \mathbf{0}_m}_{j-1\text{ times}} & 2(j-1)I_m & \underbrace{\mathbf{0}_m \; \cdots \; \mathbf{0}_m}_{n-j\text{ times}}\\
            \vdots \qquad \vdots & \vdots & \vdots \qquad \vdots\\
        \end{bmatrix} \vl. 
        \]
    Since within the matrix $C^\pm(m)$ there are only two non-zero entries of one in each row and that the entries in $\vl$ are in ``blocks'' of length $m$, we can consider the product of $C^\pm(m,n)$ with $\vl$ ``block-wise'' where the analogous ``scalar'' multiplication results in double the entries of $\vl$ in the resulting product vector $C^\pm(m,n) \vl$. Similarly, we can multiply ``block-wise'' with $D(m,n)$ and $\vl$. We now examine the cases based on the values of $1 \leq j \leq n$.
      
    \textbf{When $1 \leq j \leq i$, entries of $C^\pm(m,n) \vl$} are all  
    \begin{align*}
        \vec{0}_m^{\;\top} &= C^\pm(m) \vec{0}_m^{\;\top} \cdot i + C^\pm(m)\overrightarrow{(-1)}_m^\top \cdot (n-2i) + C^\pm(m)\overrightarrow{(n-2i)}_m^\top \\ &\qquad+ C^\pm(m) \vec{0}_m^{\;\top} \cdot (i-j) + \mathbf{0}_m \vec{0}_m^{\;\top} \cdot(j-1)\\
        &= \vec{0}_m^{\;\top} - \overrightarrow{2(n-2i)}_m^\top + \overrightarrow{2(n-2i)}_m^\top + \vec{0}_m^{\;\top} + \vec{0}_m^{\;\top}.
    \end{align*} 
    The entries of $D(m,n) \vl$ are all 
    \begin{align*} 
        \vec{0}_m^{\;\top} &= \mathbf{0}_m \vec{0}_m^{\;\top} \cdot (j-1) + 2(j-1)I_m \vec{0}_m^{\;\top} + \mathbf{0}_m \vec{0}_m^{\;\top} \cdot (i-j) \\ &\qquad + \mathbf{0}_m \overrightarrow{(-1)}_m^\top \cdot (n-2i) + \mathbf{0}_m \overrightarrow{(n-2i)}_m^\top + \mathbf{0}_m \vec{0}_m^{\;\top} \cdot (i-1)\\
        &= \vec{0}_m^{\;\top} + \vec{0}_m^{\;\top} + \vec{0}_m^{\;\top} + \vec{0}_m^{\;\top} + \vec{0}_m^{\;\top} + \vec{0}_m^{\;\top}.
    \end{align*}
    Thus, the entry in $M^\pm(m,n) \vl$ in position $j$ is $\vec{0}_m^{\;\top}$.     

    \textbf{When $i < j \leq (n-i)$ entries of $C^\pm(m,n) \vl$} are all 
    \begin{align*}
        - \overrightarrow{2(n-j-i+1)}_m^\top &= C^\pm(m)\vec{0}_m^{\;\top} \cdot i + C^\pm(m)\overrightarrow{(-1)}_m^\top \cdot (n-j-i+1) \\ & \qquad + \mathbf{0}_m \overrightarrow{(-1)}_m^\top \cdot (j-i-1) + \mathbf{0}_m \overrightarrow{(n-2i)}_m^\top + \mathbf{0}_m \vec{0}_m^{\;\top} \cdot (i-1)\\
        &= \vec{0}_m^{\;\top} - \overrightarrow{2(n-j-i+1)}_m^\top + \vec{0}_m^{\;\top}
    \end{align*}
    The entries of $D(m,n) \vl$ are all
    \begin{align*}
        - \overrightarrow{2(j-1)}_m^\top &= \mathbf{0}_m \vec{0}_m^{\;\top} \cdot i + \mathbf{0}_m \overrightarrow{(-1)}_m^\top \cdot (j-i-1) + 2(j-1)I_m \overrightarrow{(-1)}_m^\top  \\ &\qquad + \mathbf{0}_m \overrightarrow{(-1)}_m^\top \cdot (n-i-j) + \mathbf{0}_m \overrightarrow{(n-2i)}_m^\top + \mathbf{0}_m \vec{0}_m^{\;\top} \cdot (i-1)\\
        &= \vec{0}_m^{\;\top} + \vec{0}_m^{\;\top} - \overrightarrow{2(j-1)}_m^\top + \vec{0}_m^{\;\top} + \vec{0}_m^{\;\top} + \vec{0}_m^{\;\top}.
    \end{align*}
    Thus, the entry in $M^\pm(m,n) \vl$ in position $j$ is $\overrightarrow{-2(n-i)}_m^\top$.
    
    \textbf{When $j=n-i+1$ the $j^{\text{th}}$ entry of $C^\pm(m,n) \vl$} is
    \begin{align*}
       \vec{0}_m^{\;\top}  &= C^\pm(m) \vec{0}_m^{\;\top} \cdot i + \mathbf{0}_m \overrightarrow{(-1)}_m^\top \cdot (n-2i) + \mathbf{0}_m \overrightarrow{(n-2i)}_m^\top + \mathbf{0}_m \vec{0}_m^{\;\top} \cdot (i-1)\\
       &= \vec{0}_m^{\;\top} + \vec{0}_m^{\;\top} + \vec{0}_m^{\;\top} + \vec{0}_m^{\;\top}.
    \end{align*}
    The $j^{\text{th}}$ entry of $D(m,n) \vl$ is
    \begin{align*}
        \overrightarrow{2(n-i)(n-2i)}_m^\top &= \mathbf{0}_m \vec{0}_m^{\;\top} \cdot i + \mathbf{0}_m \overrightarrow{(-1)}_m^\top \cdot (n-2i) + 2(j-1)I_m \overrightarrow{(n-2i)}_m^\top \\ &\qquad + \mathbf{0}_m \vec{0}_m^{\;\top} \cdot (i-1)\\
        &= \vec{0}_m^{\;\top} + \vec{0}_m^{\;\top} - \overrightarrow{2(n-2i)(j-1)}_m^\top + \vec{0}_m^{\;\top}.
    \end{align*}
    Thus, the entry in $M^\pm(m,n) \vl$ in position $j$ is $\overrightarrow{2(n-i)(n-2i)}_m^\top$.
    
    \textbf{Finally, when $n-i+1 < j \leq n$ the entries of $C^\pm(m,n) \vl$} are all 
    \begin{align*}
        \vec{0}_m^{\;\top} &= C^\pm(m) \vec{0}_m^{\;\top} \cdot (n-j+1) + \mathbf{0}_m \vec{0}_m^{\;\top} \cdot (i+j-n-1) \\ &\qquad + \mathbf{0}_m \overrightarrow{(-1)}_m^\top \cdot (n-2i) +\mathbf{0}_m \overrightarrow{(n-2i)}_m^\top + \mathbf{0}_m \vec{0}_m^{\;\top} \cdot (i-1)\\
        &= \vec{0}_m^{\;\top} + \vec{0}_m^{\;\top} + \vec{0}_m^{\;\top} + \vec{0}_m^{\;\top} + \vec{0}_m^{\;\top}.
    \end{align*}
    The entries of $D(m,n) \vl$ are all
    \begin{align*}
        \vec{0}_m^{\;\top} &= \mathbf{0}_m \vec{0}_m^{\;\top} \cdot i + \mathbf{0}_m \overrightarrow{(-1)}_m^\top \cdot (n-2i) + \mathbf{0}_m \overrightarrow{(n-2i)}_m^\top \\ &\qquad + \mathbf{0}_m \vec{0}_m^{\;\top} \cdot (i+j-n-2) + 2(j-1)I_m \vec{0}_m^{\;\top} + \mathbf{0}_m \vec{0}_m^{\;\top} \cdot (n-j) \\
        &= \vec{0}_m^{\;\top} + \vec{0}_m^{\;\top} + \vec{0}_m^{\;\top} + \vec{0}_m^{\;\top} + \vec{0}_m^{\;\top} + \vec{0}_m^{\;\top}.
    \end{align*}
    Thus, the entry in $M^\pm(m,n) \vl$ in position $j$ is $\vec{0}_m^{\;\top}$. 
    Therefore, it follows that $M^\pm(m,n)\vl=\lambda\vl$, as desired.   

    \textbf{Case $m\equiv0\pmod4$.} Moreover, notice that we can see the matrix $M^{\pm}(m,n)$ as a block matrix with $n\times n$ blocks, each of which is one of the following $m\times m$ circulant matrices: $\textbf{0}_{m}$ or $C^\pm(m)$. Therefore, $\lambda I_{mn}-M^{\pm}(m,n)$ can be regarded as a block matrix with $n\times n$ blocks, each of them being either the $m\times m$ circulant matrix $\mathbf{0}_{m}$, $C^{\pm}(m)$, or $(\lambda-2i)I_{m}-C^{\pm}(m)$ with $0\leq i\leq n-1$ (the latter being the case of the diagonal blocks). Therefore, by Theorem~\ref{thm:Silvester} and the Leibniz formula for $\det(\lambda I_{mn}-M^{\pm}(m,n))$ yields
    \[
    \det(\lambda I_{mn}-M^{\pm}(m,n))=\det\left(\sum_{\sigma\in S_n}\sign(\sigma)M_{1,\sigma(1)}M_{2,\sigma(2)}\cdots M_{n,\sigma(n)}\right),
    \] where $M_{i,j}$ denotes the circulant matrix at position $(i,j)$ in the block matrix $\lambda I_{mn}-M^{\pm}(m,n)$. Thus,
    \begin{align*}
        \det(\lambda I_{mn}-M^{\pm}(m,n))
        &=\det\left(\sum_{\substack{\sigma\in S_n \\ \sigma\neq \operatorname{Id}_n}}\sign(\sigma)M_{1,\sigma(1)}M_{2,\sigma(2)}\cdots M_{n,\sigma(n)}+\prod_{i=0}^{n-1}(\lambda-2i)I_{mn}\right),
    \end{align*} where $\operatorname{Id}_n$ denotes the identity permutation of $S_n$. If $\lambda=2i$ for some $i\in[0,n-1]$, $\lambda I_{mn}-M^{\pm}(m,n)$, seen as a block matrix, has one row and one column where each non-zero block is $C^\pm(m)$. Therefore, Theorem~\ref{thm:Silvester} gives that above expression has the form
    \[
    \det(C^\pm(m)\mathcal{M}),
    \] for some matrix $\mathcal{M}$. Therefore, if $\lambda=2i$ for some $i\in[0,n-1]$, then $\det(\lambda I_{mn}-M^\pm(m,n))=\det(C^\pm(m))\det(\mathcal{M})$. By Lemma~\ref{lem:multipleof4}, it follows that for any $\lambda\in\{2i:i\in[0,n-1]\}$, $\det(\lambda I_{mn}-M^\pm(m,n))=0$. So in the case that $m\equiv0\pmod 4$, the spectrum of $\pmn$ contains all even integers in the set $[0,2n]$.
\end{proof}


\section{The directed case $P(m,n)$}\label{sec:directed}

Let us now consider the case $Cay(S(m,n),R^+_n)$. In this case, the graph obtained is directed if $m>2$. Note that in the case of $m=1,2$ each prefix reversal is its own inverse while when $m>2$ the inverse of $r_i^+ \in R^+_n$ is $r^-_i$ and $r^-_i \not\in R^+_n$.

\begin{thm}\label{thm:directed}
    If $m>2$, the spectrum of $P(m,n)$ contains the set $[0,n]\setminus\{\lfloor n/2\rfloor\}$
\end{thm}
\begin{proof}
    Our proof is a description of the following eigenvalue-eigenvector pairs for $M^+(m,n)$.

Since $P(m,n)$ is $n$-regular (both in and out degree for each vertex is $n$) The eigenvalue $n$ has an eigenvector that is the $mn$-dimensional vector with all 1 entries. 

In general, the eigenvalue $n-i$, with $1\leq i< \lfloor n/2\rfloor$ has a corresponding eigenvector
\[
(\underbrace{0,\ldots,0}_{mi\text{ times}},\underbrace{-1,\ldots,-1}_{m(n-2i) \text{ times}},\underbrace{n-2i,\ldots,n-2i}_{m\text{ times}},\underbrace{0,\ldots,0}_{m(i-1)\text{ times}})^\top.
\] Furthermore, if $n$ is even, say $n=2\ell$ for some positive integer $\ell$, then the eigenvalue $\ell-i$ with $1\leq i\leq  \ell$ has a corresponding eigenvector
\[
(\underbrace{0,\ldots,0}_{m(\ell-i)\text{ times}},\underbrace{-2i+1,\ldots, -2i+1}_{m\text{ times}},\underbrace{1,\ldots,1}_{m(2i-1)\text{ times}},\underbrace{0,\ldots,0}_{m(\ell-i)\text{ times}})^\top. 
\] Moreover, if $n$ is odd, say $n=2\ell+1$ for some positive integer $\ell$, then the eigenvalue $\ell+1$ has a corresponding eigenvector 
\[
(\underbrace{0,\ldots,0}_{m\ell\text{ times}},\underbrace{-1,\ldots, -1}_{m\text{ times}},\underbrace{1,\ldots,1}_{m\text{ times}},\underbrace{0,\ldots,0}_{m(\ell-1)\text{ times}})^\top.
\]
Finally, if $n$ is odd, say $n=2\ell+1$ for some positive integer $\ell$, then the eigenvalue $\ell-i$ with $1\leq i\leq  \ell$ has a corresponding eigenvector
\[
(\underbrace{0,\ldots,0}_{m(\ell-i)\text{ times}},\underbrace{-2i,\ldots, -2i}_{m\text{ times}},\underbrace{1,\ldots,1}_{2mi\text{ times}},\underbrace{0,\ldots,0}_{m(\ell-i)\text{ times}})^\top.
\] Combining all the expressions, it follows that the spectrum of $M^+(m,n)$ contains the set $[0,n]\setminus\{\lfloor n/2\rfloor\}$. Therefore, by Proposition~\ref{prop:method}, the spectrum of $P(m,n)$ contains the set $[0,n]\setminus\{\lfloor n/2\rfloor\}$.

Once again we will verify one particular case, the other cases can be approached in a similar manner. We shall verify the case where $n = 2\ell +1$ that $\lambda = \ell-i$, for $1 \leq i \leq \ell$, is an eigenvalue of $M^+(m,n)$ with corresponding eigenvector 
    
    \begin{align*}
        \vl &= (\underbrace{0,\ldots,0}_{m(\ell-i) \text{ times}},\underbrace{-2i,\ldots,-2i}_{m\text{ times}}\,\underbrace{1,\ldots,1}_{2mi\text{ times}},\underbrace{0,\ldots,0}_{m(\ell-i)\text{ times}})^\top\\
        &= (\underbrace{\vec{0}_m, \ldots, \vec{0}_m}_{\ell-i\text{ times}}, \overrightarrow{(-2i)}_m, \underbrace{\vec{1}_m, \ldots, \vec{1}_m}_{2i\text{ times}}, \underbrace{\vec{0}_m, \ldots, \vec{0}_m}_{(\ell-i)\text{ times}})^\top.
    \end{align*}
    
    To show $M^+(m,n) \vl=\lambda \vl = C^+(m,n) \vl + \frac{1}{2}D(m,n) \vl$ we argue by cases depending on the entries of $M^+(m,n) \vl $ by considering an arbitrary row $j$ of the $M^+(m,n)$, considered as a block matrix. 

        \[\begin{bmatrix}
            \vdots \qquad \vdots & \vdots \qquad \vdots\\
            \underbrace{C^+(m) \; \cdots \; C^+(m)}_{2\ell-j+2\text{ times}} & \underbrace{\mathbf{0}_m \; \cdots \; \mathbf{0}_m}_{j-1\text{ times}}\\
            \vdots \qquad \vdots & \vdots \qquad \vdots\\
        \end{bmatrix} \vl + \begin{bmatrix}
            \vdots \qquad \vdots & \vdots & \vdots \qquad \vdots\\
            \underbrace{\mathbf{0}_m \; \cdots \; \mathbf{0}_m}_{j-1\text{ times}} & (j-1)I_m & \underbrace{\mathbf{0}_m \; \cdots \; \mathbf{0}_m}_{2\ell-j+1\text{ times}}\\
            \vdots \qquad \vdots & \vdots & \vdots \qquad \vdots\\
        \end{bmatrix} \vl.
        \]
    The matrix $C^+(m)$ has only one non-zero entries of one in each row and we consider $\vl$ as a ``block vector'' with $n$ blocks of length $m$, we can consider the product of $C^+(m,n)$ with $\vl$ ``block-wise'' where the analogous ``scalar'' multiplication results in just the entries of $\vl$ in the resulting product vector $C^+(m,n) \vl$. The product of $\frac{1}{2}D(m,n)$ and $\vl$ can also be done ``block-wise.'' We now examine the cases based on the values of $1 \leq j \leq 2\ell+1 = n$.
      
    \textbf{When $1 \leq j \leq \ell-i$, entries of $C^+(m,n) \vl$} are all  
    \begin{align*}
        \vec{0}_m^{\;\top} &= C^+(m) \vec{0}_m^{\;\top} \cdot (\ell-i) + C^+(m)\overrightarrow{(-2i)}_m^\top + C^+(m)\vec{1}_m^{\;\top} \cdot (2i) \\ &\qquad + C^+(m) \vec{0}_m^{\;\top} \cdot (\ell-j-i+1) + \mathbf{0}_m + \vec{0}_m^{\;\top} \cdot(j-1)\\
        &= \vec{0}_m^{\;\top} + \overrightarrow{(-2i)}_m^\top + \overrightarrow{(2i)}_m^\top + \vec{0}_m^{\;\top}.
    \end{align*} 
    The entries of $\frac{1}{2}D(m,n) \vl$ are all 
    \begin{align*} 
        \vec{0}_m^{\;\top} &= \mathbf{0}_m \vec{0}_m^{\;\top} \cdot (j-1) + (j-1)I_m \vec{0}_m^{\;\top} + \mathbf{0}_m \vec{0}_m^{\;\top} \cdot (\ell-i-j) \\ &\qquad + \mathbf{0}_m \overrightarrow{(-2i)}_m^\top + \mathbf{0}_m \vec{1}_m^{\;\top} \cdot (2i) + \mathbf{0}_m \vec{0}_m^{\;\top} \cdot (\ell-i)\\
        &= \vec{0}_m^{\;\top} + \vec{0}_m^{\;\top} + \vec{0}_m^{\;\top} + \vec{0}_m^{\;\top} + \vec{0}_m^{\;\top} + \vec{0}_m^{\;\top}.
    \end{align*}
    Thus, the entry in $M^+(m,n) \vl$ in position $j$ is $\vec{0}_m^{\;\top}$.

    \textbf{When $j=\ell-i+1$ the $j^{\text{th}}$ entry of $C^+(m,n) \vl$} is
    \begin{align*}
       \vec{0}_m^{\;\top}  &= C^+(m) \vec{0}_m^{\;\top} \cdot (\ell-i) + C^+(m) \overrightarrow{(-2i)}_m^\top + C^+(m) \vec{1}_m^{\;\top} \cdot (2i) + \mathbf{0}_m \vec{0}_m^{\;\top} \cdot (\ell-i)\\
       &= \vec{0}_m^{\;\top} + \overrightarrow{(-2i)}_m^\top + \vec{(2i)}_m^\top + \vec{0}_m^{\;\top}.
    \end{align*}
    The $j^{\text{th}}$ entry of $\frac{1}{2}D(m,n) \vl$ is
    \begin{align*}
        \overrightarrow{(-2i)(\ell-i)}_m^\top &= \mathbf{0}_m \vec{0}_m^{\;\top} \cdot (\ell-i) + (j-1)I_m \overrightarrow{(-2i+1)}_m^\top + \mathbf{0}_m \overrightarrow{(1)}_m^\top \cdot (2i) \\ &\qquad + \mathbf{0}_m \vec{0}_m^{\;\top} \cdot (\ell-i)\\
        &= \vec{0}_m^{\;\top} + \overrightarrow{(-2i)(j-1)}_m^\top + \vec{0}_m^{\;\top} + \vec{0}_m^{\;\top}.
    \end{align*}
    Thus, the entry in $M^+(m,n) \vl$ in position $j$ is $\overrightarrow{(\ell-i)(-2i)}_m^\top$.
    
    \textbf{When $\ell - i + 1 < j \leq \ell + i + 1$ entries of $C^+(m,n) \vl$} are all 
    \begin{align*}
        \overrightarrow{(\ell-j-i+1)}_m^\top &= C^+(m)\vec{0}_m^{\;\top} \cdot (\ell-i) + C^+(m)\overrightarrow{(-2i)}_m^\top + C^+(m) \vec{1}_m^{\;\top} \cdot (\ell-j+i+1)\\ & \qquad + \mathbf{0}_m \vec{1}_m^{\;\top} \cdot (j+i-\ell-1) + \mathbf{0}_m \vec{0}_m^{\;\top} \cdot (\ell-i)\\
        &= \vec{0}_m^{\;\top} + \overrightarrow{(-2i)}_m^\top + \overrightarrow{(\ell-j+i+1)}_m^\top + \vec{0}_m^{\;\top} + \vec{0}_m^{\;\top}
    \end{align*}
    The entries of $\frac{1}{2}D(m,n) \vl$ are all
    \begin{align*}
        \overrightarrow{(j-1)}_m^\top &= \mathbf{0}_m \vec{0}_m^{\;\top} \cdot (\ell-i) + \mathbf{0}_m \overrightarrow{(-2i)}_m^\top + \mathbf{0}_m \vec{1}_m^{\;\top} \cdot (j+i-\ell-2) \\ &\qquad + (j-1)I_m \vec{1}_m^{\;\top} + \mathbf{0}_m \vec{1}_m^{\;\top} \cdot (\ell-j+i+1) + \mathbf{0}_m \vec{0}_m^{\;\top} \cdot (\ell-i)\\
        &= \vec{0}_m^{\;\top} + \vec{0}_m^{\;\top} + \vec{0}_m^{\;\top} + \overrightarrow{(j-1)}_m^\top + \vec{0}_m^{\;\top} + \vec{0}_m^{\;\top}.
    \end{align*}
    Thus, the entry in $M^+(m,n) \vl$ in position $j$ is $\overrightarrow{(\ell-i)}_m^\top$.
    
    \textbf{Finally, when $\ell+i+1 < j \leq 2\ell+1=n$ the entries of $C^+(m,n) \vl$} are all 
    \begin{align*}
        \vec{0}_m^{\;\top} &= C^+(m) \vec{0}_m^{\;\top} \cdot (2\ell-j+2) + \mathbf{0}_m \vec{0}_m^{\;\top} \cdot (j-i-\ell-2) \\ &\qquad + \mathbf{0}_m \overrightarrow{(-2i)}_m^\top + \mathbf{0}_m \vec{1}_m^{\;\top} \cdot (2i) + \mathbf{0}_m \vec{0}_m^{\;\top} \cdot (\ell-i)\\
        &= \vec{0}_m^{\;\top} + \vec{0}_m^{\;\top} + \vec{0}_m^{\;\top} + \vec{0}_m^{\;\top} + \vec{0}_m^{\;\top}.
    \end{align*}
    The entries of $\frac{1}{2}D(m,n) \vl$ are all
    \begin{align*}
        \vec{0}_m^{\;\top} &= \mathbf{0}_m \vec{0}_m^{\;\top} \cdot (\ell-i) + \mathbf{0}_m \overrightarrow{(-2i)}_m^\top \cdot (n-2i) + \mathbf{0}_m \vec{1}_m^{\;\top} \cdot (2i) \\ &\qquad + \mathbf{0}_m \vec{0}_m^{\;\top} \cdot (j-i-\ell-2) + 2(j-1)I_m \vec{0}_m^{\;\top} + \mathbf{0}_m \vec{0}_m^{\;\top} \cdot (2\ell-j+1) \\
        &= \vec{0}_m^{\;\top} + \vec{0}_m^{\;\top} + \vec{0}_m^{\;\top} + \vec{0}_m^{\;\top} + \vec{0}_m^{\;\top} + \vec{0}_m^{\;\top}.
    \end{align*}
    Thus, the entry in $M^+(m,n) \vl$ in position $j$ is $\vec{0}_m^{\;\top}$. 
    Therefore, it follows that $M^+(m,n)\vl=\lambda\vl$, as desired.
\end{proof}

\section{Final remarks}\label{sec:final}

\subsection{Multiplicities} A  question of further interest is determining the multiplicities of the eigenvalues in the spectrum. In the special case where the Cayley graph $Cay(Gr,S)$ is \emph{normal}, meaning that $S$ is closed under conjugation, then one can compute the eigenvalues and multiplicities of the adjacency matrix of $Cay(Gr,S)$ utilizing the characters of all irreducible representations of $Cay(Gr,S)$ (see, for example, Zieschang~\cite{Zieschang88}). Unfortunately, as natural as the prefix reversals are, they are not closed under conjugation. 

Computer calculations, completed to verify the above results, allow us to make the following conjecture, though. 

\begin{conj}
   If $m\equiv 0\pmod 4$, then the multiplicities of each eigenvalue of the form $2k$, with $k \in [n]\setminus \{\lfloor n/2 \rfloor$, 
   is at least 3. Moreover, the multiplicity of the eigenvalue $2\lfloor n/2\rfloor$ is at least 2. 
\end{conj}

\subsection{Spectral gap and second-largest eigenvalue}

The \emph{spectral gap} of a graph $G$ refers to either the difference between the largest eigenvalue and the second-largest eigenvalue of the adjacency matrix $A_G$ of $G$, or the smallest positive eigenvalue of the \emph{Laplacian matrix} $L_G$ of $G$ (given that the smallest eigenvalue of $L_G$ is always zero). We will consider the former definition, and denote the spectral gap of $G$ by $sp(G)$.

In the case $m=1$, the spectral gap of prefix-reversal graphs is known to be 1, see Cesi or Chung and Tobin~\cite{Cesi09, CT17}. In the case $m>1$, the spectral gap of both $P(m,n)$ and $\pmn$ is unknown. Nonetheless, in~\cite{BB24}, the authors establish that the spectral gap for the case $m=2$ is at most 1 and conjecture that $sp(P(2,n))\to1$ as $n\to\infty$. Theorem~\ref{thm:undirected} shows that $sp(\pmn)\leq 2$ and Theorem~\ref{thm:directed} shows that $sp(P(m,n))\leq 1$. Computationally, we observe that the spectral gap increases for a fixed value of $m$ as $n$ increases, which leads to the following conjecture. 

\begin{conj}\label{conj:sp}
    For fixed $m>1$, $sp(\pmn)\to2$ as $n\to\infty$ and $sp(P(m,n))\to1$ as $n\to\infty$
\end{conj}   

A closely-related question is to give an explicit formula for the second largest eigenvalue of $\pmn$ and $P(m,n)$. In the case $m=1$, the largest eigenvalue of $P_n$ is $n-1$ since $P_n$ is $(n-1)$-regular. Moreover, the work of Cesi~\cite{Cesi09} and Chung and Tobin~\cite{CT17} established that the second-largest eigenvalue of the pancake graph $P_n$ is $n-2$ (the work of Chung and Tobin is more general and addresses a family of graphs that include $P_n$). 

\subsection{Expansion ratio}

We follow the excellent survey by Hoory, Linial, and Wigderson~\cite{Hoory06}. Given a graph $G=(V,E)$ and $S\subseteq V$, the \emph{boundary} of $S$, denoted by $\partial S$, is defined as the set 
\[\partial S=\{e\in E\mid e\text{ has one end point in } S \text{ and one not in }S\}.\]

Then the \emph{expansion ratio} of $G$ is defined as
\[
h(G)=\min_{\left\{S\subseteq V: |S|\leq \frac{|V|}{2}\right\}}\frac{|\partial S|}{|S|}.
\] 

The discrete Cheeger inequality (see~\cite[Theorem 2.4]{Hoory06}) states that if $G$ is $d$-regular with adjacency matrix eigenvalues $\lambda_1(G)\geq \lambda_2(G)\geq\cdots\geq \lambda_{|V|}(G)$, then 
\[
\frac{d-\lambda_2(G)}{2}\leq h(G)\leq \sqrt{2d(d-\lambda_2(G))}.
\] Since $\lambda_2(\pmn)\geq 2(n-1)$, it follows that $h(\pmn)\leq 2\sqrt{2n}$. Notice that this  bound is independent of $m$.

Following~\cite[Definition 2.2]{Hoory06}, a sequence of $d$-regular graphs $\{G_i\}_{i\geq 1}$ of size increasing with $i$ is a \emph{family of expander graphs} if there exists $\varepsilon>0$ such that $h(G_i)\geq\varepsilon$ for all $i$. Notice that $\pmn$ is a $d_n$-regular graph with degree $d_n=n-1$ if $m=1$ and $d_n=2n$ if $m>1$. While $d_n$ increases with $n$, it is independent of $m$. Therefore, it is a reasonable question to ask if the family $\{\pmn\}_{m>1}$, with fixed $n\geq1$ and indexed by $m$, is a family of expander graphs. We conjecture that this is the case. 

\begin{conj}\label{conj:expander}
    Given $n\geq1$, the sequence of $2n$-regular graphs $\{\pmn\}_{m>1}$ indexed by $m$ is an expander family of graphs. 
\end{conj}
Indeed, if Conjecture~\ref{conj:sp} is true, then Conjecture~\ref{conj:expander} is also true. 

\bibliographystyle{acm}

\end{document}